\documentclass[12pt]{article}

\usepackage{amsmath,amsthm,amsfonts,amscd,amssymb}
\usepackage[all]{xy}

\newcommand{\M}[1]{\mathbb M_{#1}(\mathbb{C})}

\newcommand{\D}{\widehat{\Delta}}

\newcommand{\Tr}{{\rm Tr \rm}}
\newcommand{\tr}{{\rm tr \rm}}

\newcommand{\Ad}{{\rm Ad \rm}}
\newcommand{\diag}{{\rm diag \rm}}

\newtheorem{Thm}{Theorem}[section]
\newtheorem{Prop}[Thm]{Proposition}
\newtheorem{Lem}[Thm]{Lemma}

\newtheorem{Thm-f}{Théorème}[section]
\newtheorem{Prop-f}[Thm]{Proposition}
\newtheorem{Lem-f}[Thm]{Lemme}
\newtheorem{Cor-f}[Thm]{Corollaire}

\newtheorem*{Thm*}{Theorem}
\newtheorem*{Cor*}{Corollary}

\newcounter{ploum}
\newcounter{ex}[section]
\newcounter{rem}[section]
\newcounter{ass}[section]

\numberwithin{equation}{section}

\newenvironment{ass}
{\addtocounter{ass}{1}\setlength{\topsep}{1em}\par\trivlist\item{\bf Assumption
\arabic{section}.\arabic{ass}.} }{\endtrivlist}

\newenvironment{demo}[1][Proof]
{\setlength{\topsep}{1em}\par\trivlist\item{\em #1.} }
{\openbox\endtrivlist}

\begin{document}

\title{Martin Boundary Theory of some Quantum Random Walks}
\date{DMA, Ecole Normale Sup\'erieure\\
45, rue d'Ulm, 75230 Paris Cedex 05 France}
\author {Beno\^\i{}t Collins}
\maketitle

{\sc Abstract: \sc}
In this paper we define a general setting for Martin boundary theory associated
to quantum random walks, and prove a general representation theorem. 
We show that in the dual of a simply connected Lie subgroup of $U(n)$, 
the extremal Martin boundary is homeomorphic to a sphere. 
Then, we investigate restriction of quantum random walks to 
Abelian subalgebras of group algebras,
and establish a Ney-Spitzer theorem for an elementary random walk on the
fusion algebra of $SU(n)$, generalizing a previous result of Biane.
We also consider the restriction of a quantum random walk 
on $SU_q(n)$ introduced by Izumi
to two natural Abelian subalgebras, and relate the underlying Markov chains
by classical probabilistic processes. This result 
generalizes a result of Biane.

~\\
~\\

{\sc R\'esum\'e : \sc}
Dans cet article, nous d\'efinissons un cadre g\'en\'eral pour la th\'eorie
de Martin associ\'ee \`a une large classe de marches au hasard sur le dual de
groupes compacts, et \'etablissons un th\'eor\`eme de repr\'esentation int\'egrale.
Ensuite, nous montrons que dans le dual d'un sous-groupe de Lie simplement connexe
de $U(n)$, la fronti\`ere de
Martin extremale est hom\'eomorphe \`a une sph\`ere.
Nous nous concentrons alors sur la restriction de marches au hasard quantiques
\`a certaines sous-alg\`ebres Ab\'eliennes d'alg\`ebres de groupes, 
et \'etablissons un th\'eor\`eme de Ney-Spitzer pour une marche au hasard 
``de Bernoulli'' sur l'alg\`ebre de fusion de $SU(n)$.
Nous consid\'erons aussi la restriction d'une marche au hasard quantique 
introduite par Izumi \`a deux sous-alg\`ebres ab\'eliennes distinctes, et
relions les cha\^\i{}nes de Markov sous-jacentes par des proc\'ed\'es probabilistes
classiques. Ce r\'esultat g\'en\'eralise un r\'esultat de Biane.

\section{Introduction}

The classical Martin boundary theory gives a geometric
and probabilistic solution to the problem of describing positive harmonic 
functions with respect to a transient (sub)Markov operator.
This theory is well established in the framework of random walks on classical
structures. 

On the other hand, quantum probability became a self contained area in the 80's,
short after Connes developped his non commutative geometry theory, and in the
same spirit. For references, see for example
\cite{MR91k:46071}, \cite{MR91b:46061a}, \cite{MR91b:46062}, 
\cite{MR91b:46063}, \cite{MR91b:46064}, \cite{MR91b:46065}, \cite{MR91b:46066}.
Quantum probability quickly developped many 
independent and very active ramification
such as quantum stochastic processes (\cite{MR84m:82031}), 
free probability  (\cite{MR94c:46133}), and of course random matrix theory.

The problem of non-commutative harmonic analysis has also been developped for operator
algebraic and group theoretic purposes, but its developments for its 
own probabilist sake are only at its first babblings. In particular, the study 
of quantum Martin boundary was initiated by the series of papers of Biane 
(\cite{MR91k:46071}, 
\cite{MR93a:46119}, \cite{MR93b:81127}, \cite{MR94a:46091}, \cite{MR94c:46129},
\cite{MR95a:60103}) of the early nineties,
but many interesting questions raised by his papers have been left unanswered
since then. 
The purpose of this paper is to answer some of these questions, define a reasonably
general framework for quantum Martin boundary in which an integral
representation theorem works (Theorem \ref{representation}), 
and provide concrete examples. 
Note by passing that a similar representation theorem in a different
framework was obtained 
in a recent preprint of Neshveyev and Tuset \cite{math.OA/0209270}.
Amongst other applications,
we establish a purely classical Ney-spitzer theorem (Theorem \ref{martinbd})
on the set of irreducible
representation classes of the compact group $SU(n)$. The striking point of this
result is that the proof intrinsically uses results of quantum probability theory.

Theorem \ref{martinbd}
is also the starting point to an interesting counterpart of a result of
\cite{MR93a:46119}. On the group von Neumann algebra of $SU(n)$, 
the convolution operator by the normalized trace on the fundamental
representation, leaves both the center of the von Neumann algebra and the
von Neumann algebra of any maximal torus invariant. Therefore it induces two
Markov chains on discrete spaces, and they are related by an $h$-transform
(\cite{MR93a:46119}).
This result turns out to have an analogue if one replaces $SU(n)$ by
the quantum group $SU_q(n)$ of Woronowicz. In this case, classical
results about combinatorics of representation of $SU_q(n)$,
and the application of Theorem \ref{martinbd} result in 
Theorem \ref{transpro-b}.

This paper is organized as follows. 
In order to fix some notation, and for mathematical motivation, we start with
a couple of reminders. In section \ref{class-mb}, we state the Ney-Spitzer theorem
as we intend to generalize it.
Then, in section \ref{hopf-a},
we define the framework of Hopf algebras
in which we study noncommutative probability theory.
Part \ref{qmbt}
is devoted to defining a quantum Martin boundary theory and showing that any positive
harmonic element can be represented with respect to an adapted Martin kernel. 
By passing, we show that the minimal Martin boundary is isomorphic to a 
sphere under weak assumptions.
In Part \ref{nsfa}, we consider the restriction of a ``Bernoulli'' quantum walk
to the center of the Hopf algebra and establish a Ney-Spitzer like theorem.
Part \ref{fusion}
is an application of Part \ref{nsfa} to a quantum random walk on the dual of
$SU_q(n)$ introduced by Izumi (\cite{MR1916370}).

Acknowledgments: this work was completed during the author's PhD. The author
would like to thank his advisor P. Biane for many useful conversations. 
He also thanks M. Izumi for explaining him and communicating early versions of
\cite{Iz}, and R. Vergnioux for explanations about $SU_q(n)$.

\section{Reminders of (non) commutative harmonic analysis}

\subsection{Classical Martin Boundary}\label{class-mb}

Let $\mathcal{E}$ be a discrete countable state space and $P$ be a Markov
(resp. submarkov) operator defined by \index{$\mathcal{E}$}
$Pf(x)=\sum_{y\in\mathcal{E}}P(x,y)f(y)$, where $P(x,y)$ is an array of positive
real numbers assumed to satisfy $P1=1$ (resp. $P1\leq 1$).\index{$P$}
In order to avoid technical difficulties, we assume that for each $x$, every
$P(x,y)$ is zero except finitely many of them. 
Throughout the whole paper, we will identify, whenever relevant, the operator
$P$ and the kernel $P(x,y)$.
One defines inductively $P^0(x,y)=\delta_{x,y}$, $P^{n+1}(x,y)=
\sum_{z\in\mathcal{E}}P^n(x,z)P(z,y)$ and the Green kernel
$U(x,y)=\sum_{n\in\mathbb{N}}P^n(x,y)$.
We make the usual assumptions of irreduciblity and 
transience: \index{$U$}\index{noyau!de Poisson}

\begin{ass}\label{martin}
For all $x,y\in\mathcal{E}$, one has $0<U(x,y)<\infty$. 
\end{ass}

A function $f$ is said to be {\it harmonic\it} with respect to $P$
if $Pf=f$.
The {\it Martin kernel \it} with a base point $x_0$ is 
defined as \index{noyau!de Martin}
\begin{equation*}\index{$K$}
k(x,y)=U(x,y)/U(x_0,y)
\end{equation*}
Harnack inequalities imply 
that for all $x\in\mathcal{E}$, the function $k(x,.)$ is bounded.
The {\it Martin compactification \it} $MS$ of $\mathcal{E}$ is defined as the smallest compact
subspace in which $\mathcal{E}$ can be continuously and densely embedded and such that every
function $k(x,.)$ can be uniquely extended by continuity to $MS$. 
Let $MB$ be the boundary of $\mathcal{E}$ in $MS$. \index{$MB$}\index{$MS$}
\index{Martin!espace de}\index{Martin!frontière de}
A positive harmonic function $f$ is said to be {\it minimal\it} if 
any harmonic function $g$ satisfying $0\leq g\leq f$ is a multiple of $f$.
There exists a measurable subset $MB^{min}$ of $MB$ such that $x\in MB^{min}$
if and only if $k(.,x)$ is a positive minimal harmonic function.

\begin{Thm} \label{frontieremartin}\index{$MB^{ex}$}
For every positive harmonic
function $f$ satisfying $f(e)=1$, there exists a unique probability measure
$\mu_f$ on $MB$ such that $\mu_f(MB^{min})=1$ and 
for all $x\in\mathcal{E}$, $f(x)=\int_{\xi\in MB^{min}}
k(x,\xi)d\mu_f(\xi)$.
\end{Thm}

For the proof of the above theorem, see \cite{MR86a:60097}
or \cite{MR53:11748}.
In view of this, it is natural to try to compute explicit examples. This turns out
to be a difficult task, and one remarkable example of such a computation is done in
\cite{MR33:3354}. 

Let $\mathcal{E}$ be the state space $\mathbb{Z}^d$, for some $d\geq 2$, and
$\mu$ be a finitely supported measure (the hypothesis of finite 
support can be considerably weakened, but we do not enter into such technical 
considerations) whose mean on $\mathbb{R}^d$ is
different from zero. We identify canonically $\mathbb{Z}^d$ with a lattice of 
the Euclidean space $\mathbb{R}^d$ with its scalar product $\langle .,. \rangle$. 
Then, the set 
\begin{equation*}
E=\{x, \int_{\mathbb{Z}^d} \exp \langle x,X \rangle d\mu (X)=1\}
\end{equation*} 
is a $C^{\infty}$ submanifold of $\mathbb{R}^d$. It is 
diffeomorphic to the 
sphere $S^{d-1}$. It is a consequence on a theorem of Choquet and Deny 
(see \cite{MR22:9808})
that any positive function on $\mathbb{Z}^d$ harmonic
with respect to the operator of convolution by
$\mu$ (we call its operator $P_{\mu}$) admits an unique integral representation as a 
linear combination of functions $y\rightarrow e^{(x,y)}$, $x\in E$.\index{$P$}

Let $\beta$ be a continuous increasing bijection from $[0,\infty ]$ to
$[0,1]$.
Let $\Pi$ be the map from $\mathbb{Z}^d$ to the canonical unit ball $B(0,1)$
of $\mathbb{R}^d$ given by $\Pi (x)=\beta (||x||) x/||x||$. 
This map is a topological injection, and
the metrics it inherits gives rise to a compactification of $\mathbb{R}^d$ by $S^{d-1}$.
The map $\Pi$ extends continuously
to $\mathbb{Z}^d\cup S^{d-1}$ with value in $B(0,1)$.
The map $u$ sending an element of $E$ to its normed outer normal vector is a $C^{\infty}$ 
diffeomorphism from $E$ to $S^{d-1}$.

\begin{Thm}[NS, \cite{MR33:3354}]
The Martin compactification of $\mathbb{Z}^d$ is $\mathbb{Z}^d\cup S^{d-1}$.
The correspondence between $S^{d-1}$ and $E$ arising from
the Martin kernel is given by the map $u$.
The extremal Martin boundary $MB^{min}$ thus coincides with $MB$.  
\end{Thm} 

Throughout the paper, we define
\begin{eqnarray}
\mathcal{H}^+_{\mathcal{E},P}=\{f\in\mathbb{C}^{\mathcal{E}}, f\geq 0, Pf=f, f(e)=1\}\\
\mathcal{H}_{\mathcal{E},P}=\{f\in L^{\infty}(\mathcal{E}), Pf=f \}
\end{eqnarray}
\index{$\mathcal{H}^+_{\mathcal{E},P}$}\index{$\mathcal{H}_{\mathcal{E},P}$}

The set $\mathcal{H}^+_{\mathcal{E},P}$ is convex compact for the topology of
pointwise convergence. Therefore it admits extremal points.
By definition, let $\mathcal{H}^{+,ex}_{\mathcal{E},P}$ be this subset.
By the Krein-Milman theorem, the closure of the convex hull of these extremal points
is exactly $\mathcal{H}^+_{\mathcal{E},P}$.\index{$\mathcal{H}^{+,ex}_{\mathcal{E},P}$}

Thus, Theorem \ref{frontieremartin} identifies the extremal points of the above
set with $MB^{min}$. 
We shall say in this paper that we obtain a {\it Ney-Spitzer like\it} theorem when
we give a compactification $MS=\mathcal{E}\cup MB$ 
of a state space, describe the extremal positive
harmonic functions, describe a subset $MB^{min}$ of $MB$ and a bijection between
$\mathcal{H}^{+,ex}_{\mathcal{E},P}$ and $MB^{min}$.

  \subsection{Hopf algebras}\label{hopf-a}

We fix some classical notations of operator algebra theory, and 
remind some elementary definitions and results of Hopf algebra theory.

Let $G$ be a topological compact group, and $d\mu$ its left and right invariant
probability Haar measure.
Let $L^2(G)$ be the $L^2$ space associated to this Haar measure.
For $g\in G$, the unitary operator
$\lambda_g\in B(L^2(G))$ is defined by
$\lambda_g : f\rightarrow
(x\rightarrow f(g^{-1}x))$. The mapping $g\rightarrow\lambda_g$ is 
continuous for the strong operator topology in $B(L^2(G))$.
The vector space $Vect (\lambda_g , g\in G)$ is
a $*$ -subalgebra of $B(L^2(G))$. 
Let $M(G)$ be the von Neumann algebra of $G$, i.e. the bicommutant of
$Vect (\lambda_g , g\in G)$ in $B(L^2(G))$.\index{$M(G)$}

The set of equivalence classes of irreducible finite dimensional
unitary representations of $G$ is denoted by $\Gamma$. 
For $x\in\Gamma$, $d_x$ is the corresponding dimension.
By Peter-Weyl's theorem, we have the following isomorphism of von 
Neumann algebras:

\begin{equation*}
M(G) \cong \oplus_{x\in\Gamma} \M{d_x}
\end{equation*}

The counit is the map $M(G)\rightarrow \mathbb{C}$, \index{$\varepsilon$}
defined as the continous linear expansion of
the map $\varepsilon (\lambda_g)=1$. 

In the same way, the antipode \index{$S$} 
is the continous map $M(G)\rightarrow M(G)$
such that $S(\lambda_g ) =\lambda_{g^{-1}}$, and the coproduct is the map
$M(G)\rightarrow M(G\times G)\cong M(G)\otimes M(G)$ such that
$\D (\lambda_g) =\lambda_g\otimes \lambda_g$. 

The quadruple $(M(G),\varepsilon , \D, S)$ is called the Hopf-von Neumann algebra
of the group $G$. \index{counité}\index{algèbre!de Hopf}\index{coproduit}

For a von Neumann algebra $A$, we shall call $\widehat{A}$ the set of elements
affiliated to $A$.\index{$\widehat{A}$}
By the Peter-Weyl theorem, $\widehat{M(G)}$ is a $*$ -algebra endowed
with a natural pointwise convergence topology. As a topological $*$ -algebra, it
is isomorphic to $\prod_{x\in\Gamma}\M{d_x}$. 

One can see (see \cite{MR94a:46091}), that
$\D$, 
$\varepsilon$ and $S$ are also continuous for the topology of pointwise convergence, 
so that there is a unique way of extending them from $\widehat{M(G)}$
to $\widehat{M(G)\otimes M(G)}$ (resp., $\widehat{M(G)}$, $\mathbb{C}$).
According to Effros and Ruan (see \cite{MR95j:46089}; also see \cite{MR94a:46091}) we call 
this structure a {\it topological $*$ -Hopf algebra\it}.

Let $\overset{o}{M(G)}$ be the set of finite rank operators in $M(G)$.
A linear form $\nu :\widehat{M(G)}\rightarrow\mathbb{C}$ is said to
be {\it finitely supported\it} iff it is continuous with respect to the
pointwise convergence. Equivalently, there exists a faithful weight $\tau$ on
$M(G)$ and an element 
$A\in\overset{o}{M(G)}$, such that for all $B\in\widehat{M(G)}$, one has
$\mu (B)=\tau(AB)$. We call $(\widehat{M(G)})_*$ the vector space
of finitely supported linear forms.

To summarize, we will be dealing with the following inclusions of algebras:
\begin{equation*}
\overset{o}{M(G)}\subset M(G)\subset\widehat{M(G)}
\end{equation*}
The first one is not a Hopf algebra, but the latter two are. 

  \subsection{Random walks and harmonic analysis}

We use the framework of Hopf algebras in order to define quantum random walks. Several groups
of mathematicians have already inspected axiomatics (see \cite{MR84m:82031})
and their properties (see for example \cite{Iz}, \cite{MR2000b:81083}, \cite{MR96k:60004}).

For $l\in\Gamma$, let $1_l$ be \index{$1_l$}
the minimal central idempotent in $M(G)$ associated to the irreducible representation $l\in\Gamma$ in
$M(G)$. For $\nu$ and $\mu$ two states on $M(G)$, we define their
convolution $\mu *\nu$ by the equation
$\mu *\nu (f)=(\mu\otimes\nu)\D (f)$ for each $f\in M(G)$.
$\mu *\nu$ is a state and we can define inductively 
$\nu^{*n}$ to be $\varepsilon$ if $n=0$ and $\nu *\nu^{*n-1}$ else.
We define the operator $P_{\mu}$ on $M(G)$ by 

\begin{equation*}
P_{\mu}(f)=(id\otimes \mu)\D (f)
\end{equation*}
and its iterates inductively by $P_{\mu}^n=id$ if $n=0$ 
and $P_{\mu}^n=P_{\mu}\circ P_{\mu}^{n-1}$ else.
One has also $P_{\mu}^n=(id\otimes\mu^{*n})\D$.\index{$P_{\mu}$}

The operator $P_{\mu}$ is the evolution operator associated to a
quantum random walk on the dual of $G$ in the sense of \cite{MR84m:82031}.
It is a completely positive operator on $M(G)$. 
If $\mu (1_l)=0$ for any but finitely many $l$'s,
$P_{\mu}$ extends to a positive continuous operator on
$\widehat{M(G)}$.

An element $f$ in $\widehat{M(G)}$ is said to be 
{\it harmonic \it} with respect to $P_{\mu}$ iff \index{harmonique!\'el\'ement}
\begin{equation}
P_{\mu}f=f
\end{equation}

Biane showed in  \cite{MR93b:81127}, that 
$\mathcal{H}^{+,ex}_{P_{\mu}}$ is the set 
\begin{equation}
E=\{f\in \widehat{M(G)}, f\geq 0, \D f=f\otimes f, \mu(f)=1\}
\end{equation}

We shall say that an element $f\in \widehat{M(G)}$ such that
$\D f=f\otimes f$ is an {\it exponential\it}, and call
$Ex(G)$ be the set of exponentials.

\section{Quantum Martin boundary theory}\label{qmbt}

  \subsection{Representation of positive harmonic elements}\label{sect-repres}

In this section, we define a Martin compactification and a Martin kernel,
and show that every positive harmonic element can be represented by a state
on the Martin boundary.
For a completely positive continuous operator $Q$ from $M(G)$ into itself
and $\nu$ a weight, $\nu Q$ is again a weight defined by $\nu Q f=\nu (Q(f))$.
We write $\nu\leq \mu$ iff for any positive $f$, $\nu (f)\leq \mu (f)$.
We need the following assumption:

\begin{ass}\label{martinbis}
The weight $\mu$ is such that $\mu (1)=q\in]0,1[$ and that
$\varepsilon U=\sum_{n\geq 0} \mu^{*n}$ is faithful. 
There exists a positive $A\in Ex(G)$ and a tracial operator 
$\widetilde{\mu}$ such that for any $f\in\widehat{M(G)}$,
$\mu (f)=\widetilde{\mu}(Af)$.
\end{ass}

Let $U=\sum_{n\geq 0}P_{\mu}^n$. This operator has operator
norm less than $(1-q)^{-1}$ and 
is the quantum analogue of the {\it Green kernel\it}.\index{$U$}\index{noyau!de Poisson}
\index{$\overset{o}{M(G)}$}
Let $\overset{o}{M(G)}$ be the subalgebra of $M(G)$, consisting
of finite rank operators. Note that this is not a Hopf subalgebra.

We define the {\it Martin Kernel \it} to be the linear map 
\begin{equation}\index{noyau!de Martin}
K: (\widehat{M(G)})_*\rightarrow\widehat{M(G)}
\end{equation}
such that for any $\nu\in (\widehat{M(G)})_*$, $K_{\nu}$ 
satisfies for all $f\in\overset{o}{M(G)}$,
\begin{equation}\label{defmartin}
\nu U(f)= \varepsilon U(A^{-1/2}K_{\nu}A^{1/2}f)
\end{equation}

$K_{\nu}$ is well defined because the weight $\varepsilon U$ is faithful.
A definition equivalent to this one has already appeared in P. Biane's papers
(see \cite{MR95a:60103}) in the context of $SU(2)$ for a tracial weight.

\begin{Lem}\label{thpot}
Let $g\in M(G)$. Then one has $g=U(g-P(g))$
\end{Lem}

\begin{demo}
It is enough to remark that $Uf$ is defined for any $f\in M(G)$, that
under the Assumption \ref{martinbis} it is the norm limit of $\sum_{k=0}^n P^k f$, 
and that 
$\sum_{k=0}^n P_{\mu}^k (f-P_{\mu}f)=f-P_{\mu}^{n+1}f$. But $P_{\mu}^{n+1}f$ has norm tending towards
zero as $n$ tends to infinity.
\end{demo}

\begin{Prop} \label{propK}
\begin{itemize}
\item $K$ is positive and its image is contained in $M(G)$. 
\item The norm closure of ${\rm span \rm}\{K_{\nu},\nu\in (\widehat{M(G)})_*\}$ 
contains the $C^*$ -algebra $\mathcal{K}$ 
of compact operators. 
\end{itemize}
\end{Prop}

\begin{demo}
The weight $\varepsilon U$ is faithful. Furthermore, it is
tracial because it is invariant under the adjoint action.
It is known that $\phi : \M{n}\rightarrow\mathbb{C}$ is a weight iff 
there exists a positive matrix $B$ such that for all $A$, $\phi(A)=\Tr(AB)$.
Therefore, if $\nu$ is a weight then $K_{\nu}$ is positive.

For $\nu$ a finitely supported weight, there exists by Assumption
\ref{martinbis} an integer $n$ and
a constant $\alpha$ such that $\nu\leq \alpha\sum_{i=0}^n\varepsilon P^i$.
By positivity of $P$ this implies that $\nu U\leq (n+1) \alpha \varepsilon U$.
This implies that $K_{\nu}$ is bounded.

If $K_{\nu}=0$, then by faithfulness of $\varepsilon U$, one has
$\nu U=0$, and by Lemma
\ref{thpot}, $\nu =0$, which proves the ``into''.

For the second statement, it is enough to show that for every positive
finite dimensional operator $f\in\widehat{M(G)}$,
there exists $\nu$ such that $K_{\nu}=f$. 

Let $\nu $ be the linear form such that for any $g\in M(G)$,
$\nu (g)=\varepsilon U (f(g-P_{\mu}g))$. Then $\nu$ is finitely supported,
so that $K_{\nu}$ is well defined;
and one can check using Lemma \ref{thpot} that it satisfies $K_{\nu}=f$.
\end{demo}

We define the {\it Martin space \it} $MS$ 
to be the $C^*$-algebra
\begin{equation*}
MS=C^*(K_{\nu}, \nu\in (\widehat{M(G)})_*)
\end{equation*}
By Proposition \ref{propK}, $\mathcal{K}$ is an ideal of $MS$. 
Let the {\it Martin boundary \it}
be the $C^*$ -algebra $MB=MS/\mathcal{K}$.
The Martin compactification of the dual $\widehat{G}$ of $G$ is then
defined to be the following exact sequence:
\begin{equation}
0\rightarrow\mathcal{K}\rightarrow MS \rightarrow MB \rightarrow 0
\end{equation}

We are now able to prove the following representation theorem:

\begin{Thm} \label{representation}
\begin{itemize}
\item For each positive harmonic element $h$ in $\widehat{M(G)}$ there exists
a state $\phi_h$ on $MB$ such that for every finitely supported linear
form $\nu$ on $M(G)$, one has $\nu (h)=\phi_h (K_{\nu})$
\item This representation is unique if 
${\rm span \rm}\{K_{\nu},\nu\in (\widehat{M(G)})_*\}$ is dense in $MS$.
\end{itemize}
\end{Thm}

\begin{demo}
If $\phi$ is a state on $MB$ then the element 
$h_{\phi}\in\widehat{M(G)}$ defined by
$\nu (h_{\phi})=\phi (K_{\nu})$ for all $\nu\in(\widehat{M(G)})_*$,
is positive. Let us first show that it is harmonic.

We need to show that $P_{\mu}h_{\phi}=h_{\phi}$, or equivalently,
that for all $\nu\in(\widehat{M(G)})_*$,
$\nu (P_{\mu}h_{\phi})=\nu (h_{\phi})$.
But $\nu Ph_{\phi}=(\mu *\nu )h_{\phi} = \phi (K_{\mu *\nu})$
and $\nu (h_{\phi})=\phi (K_{\nu})$, therefore it is enough
to show that $\phi (K_{\mu *\nu})=\phi (K_{\nu})$.

We have, for all $f\in \widehat{M(G)}$
\begin{equation*}
\begin{split}
\varepsilon U((K_{\mu *\nu}-K_{\nu})f)=\\
(\mu *\nu-\nu)Uf=\nu f
\end{split}
\end{equation*}
(second equality arises from Lemma \ref{thpot}). But since
$\nu\in(\widehat{M(G)})_*$, this implies that $K_{\mu *\nu}-K_{\nu}$
has finite rank, thus is compact. Therefore
$\phi (K_{\mu *\nu}-K_{\nu})=0$

Furthermore, $h_{\phi}$ satisfies $\varepsilon h_{\phi}=1$ and 
one sees that the linear map
\begin{equation*}
\Xi : (MS)^* \rightarrow \widehat{M(G)}
\end{equation*}
that maps $\phi$ to $h_{\phi}$ is continuous for the pointwise convergence
topology.

We will now show that for every extremal harmonic element $h$, there exists a 
weight $\phi$ on $MS$ such that for any $\nu\in (\widehat{M(G)})_*$, one has
\begin{equation*}
\nu (h)=\phi (K_{\nu})
\end{equation*}
Since $\Xi$ is linear, any convex combination of extremal harmonic elements
can be represented. Furthermore, $\Xi$ is continuous, therefore, any harmonic
element that can be approximated in the pointwise convergence topology 
by a convex combination of extremal harmonic elements can be represented,
therefore, any element can be represented.

Let $h$ be such a minimal harmonic element. 
By Biane's theorem (\cite{MR94a:46091}),
it satisfies $\D h=h\otimes h$ and $\mu (h)=1$. Let $\widehat{M(h)}$ be the 
closure in $\widehat{M(G)}$ of the algebra generated by $h$. This is obviously
a topological Hopf $*$-subalgebra of $\widehat{M(G)}$.
This algebra is commutative, cocommutative and closed, therefore one can show 
directly that the operations $inf$ and $sup$ are well defined inside this algebra.

In $\widehat{M(h)}$, consider
a sequence $h_k$ of positive elements such that $U h_k$ tends and increases
towards  $h$ as $k$ goes towards infinity. The existence of such
a sequence is a consequence of standard probabilistic considerations, but we justify
it nonetheless. 

Let $f_k=inf (k Id, h)$, where the infimum is taken
on the commutative affiliated algebra $\widehat{M(h)}$. 
$f_k$ is bounded and satisfies $P_{\mu}(f_k)\leq f_k$. Let $h_k=f_k-Pf_k$.
This element is positive and it can not be zero because $h$ is extremal and
non bounded, thus $f_k$ would have to be a multiple of $h$, which would result in
$h=0$. By Lemma \ref{thpot} 
this implies that $f_k=Ug_k$. Last, it is obvious that $f_k$ tends towards $h$
in the pointwise convergence topology, as $k\rightarrow\infty$.

Consider $\phi_k=\varepsilon U(A^{-1/2}h_kA^{1/2} \cdot)$. 
It is a state on the norm closed operator system generated by $K_{\mu}$, and 
it satisfies $\phi_k (1)=\varepsilon Uh_k \leq 1$. By a classical result 
(see \cite{MR30:1404}, p. 50, lemme 2.10.1), it extends to a state on $MS$.

Furthermore, $\phi_k (K_{\nu})=\nu (Uh_k)$ tends towards $\nu (h)$ as $k$ tends towards
infinity. This proves that
$\phi_k$ converges weakly towards a state $\phi$ on $MS$.
This state vanishes on compact operators, so is actually a state of $MB$.
\end{demo}

  \subsection{Topological structure of the boundary}

In this section, we assume that $G$ is a compact simply connected 
Lie subgroup of $\mathbb{U}_n(\mathbb{C})$ with Lie algebra $\mathfrak{g}$.
Let $\mathfrak{g}_{\mathbb{C}}$ be the complexified Lie algebra, 
$G_{\mathbb{C}}$ be the complexified Lie group, and $(\rho, V)$ be
the fundamental representation of $G$.

The left regular representation yields an identification of
$\mathfrak{g}_{\mathbb{C}}$ with a Lie subalgebra of $\widehat{M(G)}$.
If $f\in \mathfrak{g}_{\mathbb{C}}$, then $\D f= f\otimes 1 + 1\otimes f$.
The map $EXP : \widehat{M(G)}\rightarrow \widehat{M(G)}$
defined by the usual series is such that
for any $f\in \mathfrak{g}_{\mathbb{C}}$,
$\D EXP f=EXP f \otimes EXP f$.

By a slight modification of a
result of Biane (\cite{MR93b:81127}, Proposition 11 and 
Lemme 12), if any irreducible 
representation of $G$ is contained in some tensor power of
$V$, then the set of non-zero exponentials in $\widehat{M(G)}$ is 
the group generated by $EXP \mathfrak{g}_{\mathbb{C}}$.
It is exactly $G_{\mathbb{C}}$, and the restriction 
of $Ex(G)$ to $End (V)$ is a group isomorphism between 
$Ex(G)$ and $G_{\mathbb{C}}$. 
An explicit isomorphism is obtained by restricting $Ex(G)$ to
the fundamental representation. We call 
\begin{equation}\label{def-i}
i : G_{\mathbb{C}} \rightarrow Ex(G)
\end{equation}
the converse of this
isomorphism.

The following theorem answers a question raised by Biane about the topology
of the boundary.

\begin{Thm} \label{sphere}
Let $\mu$ be a weight on $\widehat{M(G)}$ satisfying $\mu (1)=q\in]0,1[$ and
assumption \ref{martinbis}. 

Then, the set $H^{+,ex}_{P_{\mu}}$ of extremal harmonic elements is diffeomorphic 
to the sphere
$S^{k-1}$, where $k$ is the dimension of the Lie algebra $\mathfrak{g}$.
\end{Thm}

\begin{demo}
Let $\mathfrak{g}_{sa}$ be the real vector subspace of $\widehat{M(G)}$ 
of self adjoint elements of $\mathfrak{g}_{\mathbb{C}}$ in $\widehat{M(G)}$.
Let $x\in \mathfrak{g}_{sa}$
be non-zero, and $f_x$ the map
$\mathbb{R}\rightarrow\mathbb{R}$ given by $f_x(t)=\nu (EXP tx)$.
This map is always positive. Since $\tr x=0$ and $x$ is Hermitian,
it has one negative eigenvalue and one positive eigenvalue. Therefore
$\lim_{\pm\infty}f_x=\infty$. Besides by definition, $f_x(0)=q<1$.
The function $f_x$ admits the second derivate $\nu (x^2 EXP tx)$ at $t$, therefore
it is positive. Thus, the function $f_x$ is convex; therefore there exists only
two real numbers $t_x^+$ (resp. $t_x^-$) satisfying $t_x^+>0$ and $f_x(t_x^+)=1$
 (resp. $t_x^-<0$ and $f_x(t_x^-=1)$).
But the map $EXP$ is a diffeomorphism from $\mathfrak{g}_{sa}$ onto 
$E=\{x\in\widehat{M(G)}, \D x=x\otimes x \}$.
Therefore the inverse image of $\mathcal{H}^{+,ex}_{P}$ under $\exp$ is a closed
star-like subset around $0$, therefore it is homeomorphic to the sphere $S^{k-1}$.
\end{demo}

\section{A Ney-Spitzer theorem for a random walk on a Weyl chamber}\label{nsfa}

It would be interesting and seems challenging to obtain nice generalizations
of the result of \cite{MR95a:60103} in the framework developped above. 
We are not able to perform fully such computations.
Yet, it is possible to obtain a Ney-Spitzer like theorem if one restrict
a tracial quantum random walk on $SU(n)$ to the center of its Hopf algebra.
In this section, we establish a Ney-Spitzer like theorem for the most elementary
quantum random walk, improving previous results of \cite{MR93a:46119}.

\subsection{Main result}\label{statement}

In the Euclidean space $\mathbb{R}^n$, $n\geq 3$ with canonical basis 
$(\tilde{e_i})_{i=1}^n$, we consider the lattice $L$ spanned by\index{$L$}
$e_i=\tilde{e_i}-(\tilde{e_1}+\ldots +\tilde{e_n})/n$. 
There is a unique way to write $x\in L$ under the form
$x=\sum_{i=1}^nx_i e_i$ such that every $x_i\in\mathbb{N}$ and one at least is zero.
We call $(x_1,\ldots ,x_n)$ the {\it coordinates \it} of $x$ and $\sum x_i = |x|$ the 
{\it length \it} of $x$. \index{$|x|$}
Let\index{$\overset{o}{W}$}
\begin{equation}\index{$W$}
\begin{split}
W=\{ x\in L , x_1\geq x_2\geq\ldots\geq x_n=0\}\\
\overset{o}{W}=\{x\in L , x_1> x_2>\ldots> x_n=0\}\\
\end{split}
\end{equation}

The lattice $L$ is compactified by the sphere $S^{n-2}$ in the following sense: 
a sequence $x^d$ of $L$ tends towards $y\in S^{n-2}$ iff its Euclidean norm $||x_n||$
tends to infinity and $x^d/||x^d ||\rightarrow y$. 

Consider the simplex
\begin{equation}\label{simplexe1}
\Sigma=\{y'=(y_1',\ldots ,y_n') \, : \, 
y_1'\geq\ldots\geq y_n'=0, \sum_iy_i'=1\}\subset \mathbb{R}^n
\end{equation}

We embed it into $S^{n-2}$ by the map $y'\rightarrow y'/||y'||$.
By doing so, the above compactification induces a compactification of 
$\overset{o}{W}$ by $\Sigma$.

Namely, let $y^d=(y^d_1>y^d_2>\ldots >y^d_n=0)_{d\geq 0}$ be a series of elements 
of $\overset{o}{W}$. Then it converges iff $y_d$ is constant
for $d$ large enough or if $|y_d|$ tends towards infinity and 
for all $i$, $y^d_i/|y^d|$ admits a limit $y_i'$. 

Let $0<\alpha <1$ be a real number. Consider the measure
\begin{equation}
\mu=\sum_{i=1}^n \frac{\alpha}{n} \delta_{e_i}
\end{equation}
where $\delta_x$ is the Dirac mass at $x$. 

The Martin theory with respect to $P_{\mu}$ is completely understood.

Let us now define a new random walk on the state space $\overset{o}{W}_n$.
The vector $\rho=(n-1,n-2,\ldots, 0)$ is such that  
$W_n=\rho+\overset{o}{W}_n$. 
In the sequel, we abbreviate $P_{\mu}$ by $P$.
Our random walk in $\overset{o}{W}_n$ is obtained from
$P$ and conditioned not to hit $\partial W_n=W_n - \overset{o}{W}_n$. 
Call $\overset{o}{P}$ its transition kernel. For $x\in\overset{o}{W}_n$,
its transition kernel satisfies
$\overset{o}{P}(x,y)=\alpha /n$ iff $y-x=e_i$ and $y\in\overset{o}{W}_n$ and $0$ otherwise.
The main result of this section is:
\begin{Thm}\label{martinbd}
The Martin boundary associated to $\overset{o}{P}$, the state space
$\overset{o}{W}$ and the vector $e$, is homeomorphic to $\Sigma$. Furthermore, 
$MB=MB^{min}$.
\end{Thm}
The remainder of this section is devoted to proving this theorem.
In part \ref{mb}, we compute $MB$. In part \ref{mex}, we consider
$\overset{o}{W}$ as a canonical basis for the fusion algebra of
$SU(n)$, and consider the Markov operator $\widetilde{P}$ \index{$\widetilde{P}$}
obtained by considering the convolution by the normalized fundamental
representation. Thanks to a result of Biane, its abstract Martin boundary
can be computed and identified with that of $\overset{o}{P}$.
We use that to show that $MB=MB^{min}$.

\subsection{Asymptotics of the Martin kernel}\label{mb}

The translation invariance of $P_{\mu}$ implies that we can define for any $l$ 
the one parameter functions 
\begin{equation*}
P_{\mu}^l(y-x)=P_{\mu}^l(x,y)
\end{equation*}
From now on, we take the convention that $x!^{-1}=0$ if $x<0$.

\begin{Lem}
With the notations of section \ref{statement}, we have
\begin{itemize}
\item
\begin{equation*}
P_{\mu}^l(y)=
\begin{cases}
\frac{(|y|+kn)!}{\prod_{i=1}^n(y_i+k)!}(q/n)^l \,\,{\rm if \rm} \,\,\, l=kn+|y|\\
0 \,\, {\rm otherwise \rm}
\end{cases}
\end{equation*}
\item
\begin{equation*}
\overset{o}{P}^l(x,y)=(|y-x|+kn)!(q/n)^{|y-x|+kn}\det(y_i-x_j+k)!^{-1}
\end{equation*}
\end{itemize}
\end{Lem}

\begin{demo}
The first point is elementary combinatorics. For the second one,
if $x=(x_1,\ldots ,x_n)$, then for $\sigma$ a permutation of $[1,n]$,
let $x_{\sigma}=(x_{\sigma (1)},\ldots ,x_{\sigma (n)})$. 
We have, by the reflexion principle,
\begin{equation*}
\overset{o}{P}^l(x,y)=\sum_{\sigma\in\mathcal{S}_n}P_{\mu}^l(x_{\sigma},y)\varepsilon(\sigma )
\end{equation*}
where $\varepsilon (\sigma )$ is the signature of the permutation $\sigma$, and
the result follows.
\end{demo}

Let $Y=(Y_1,\ldots ,Y_n)$ be an $n$-tuple of formal variables.
Recall that the Vandermonde determinant is the polynomial \index{Vandermonde}
$V (Y)=\prod_{1\leq i<j\leq n}(Y_j-Y_i)$. Let $x\in W$.
The function
\begin{equation}\index{$s_{\lambda,d}$}
s_x(Y_i)=\frac{\det (Y_i^{x_j+n-j})}{V(Y_i)}
\end{equation}
is a symmetric polynomial in $Y$ homogeneous of degree $|x|-|e|$. 
It is known as the {\it Schur polynomial \it}
(see \cite{MR99f:05119}).
It is classical (see \cite{MR93a:20069}) that $W$ is in one to one
correspondance with
the set of classes of irreducible representations of $SU(n)$ up to isomorphism.
By Weyl character formula, 
$s_x$ is known to be the character of the irreducible representation
associated to $x$ evaluated on $\diag (Y_1,\ldots ,Y_n)$.

\begin{Lem}\label{lem-schur}
Let
\begin{equation*}
f(x,y,l)=\frac{(|y-x|+kn)!}{(kn+|x|+|y-x|)!}
\end{equation*}
and
\begin{equation*}
c_{x,y}=|x|+|y-x|-|y|
\end{equation*}
Note that by the triangle inequality, $0\leq c_{x,y} \leq |x|$. One has

\begin{itemize}
\item
For each $x,y\in\overset{o}{W}_n$ and $l\in\mathbb{N}$ such that $l=|y-x|+kn$ for some
positive integer $k$, 
\begin{equation*}
\overset{o}{P}^l(x,y)=(q/n)^{-|x|}V(y_i)f(x,y,l)p^{l+c_{x,y}}(y)s_{x-\rho}(y_i+k)(1+o(1))
\end{equation*}
where the symbol $o(1)$ has to be understood pointwise 
in $x$, as $|y|\rightarrow\infty$, uniformly in $k\geq 0$.
\item
$f$ is equivalent to $(|y|+kn)^{-|x|}$ as $|y|\rightarrow\infty$, uniformly
in $k\geq 0$.
Furthermore $f(x,y,l)s_{x-\rho}(y_i+k)(q/n)^{-|x|}$ is bounded independently 
on $|y|$ and $k$.
\end{itemize}
\end{Lem}

\begin{demo}
Since
\begin{equation*}
\overset{o}{P}^l(x,y)=(|y-x|+kn)!(q/n)^{|y-x|+kn}\det(y_i-x_j+k)!^{-1}
\end{equation*}
we have 
\begin{equation*}
\overset{o}{P}^l(x,y)/p^{l+c_{x,y}}(y)=(q/n)^{-|x|}f(x,y,l)\det 
(\frac{(y_i+k)!}{(y_i-x_j+k)!})
\end{equation*}
The expression  $(y_i+k)!/(y_i-x_j+k)! $ is a polynomial
in the variable $y_i+k$ whose leading term is $(y_i+k)^{x_j}$.
By multilinearity of the determinant and the definition of Schur polynomials 
this implies that 
\begin{equation*}
\det (\frac{(y_i+k)!}{(y_i-x_j+k)!})=V(y)s_{x-\rho}(y_i+k)(1+o(1))
\end{equation*}
The element $s_{x-\rho}(y_i+k)$ is positive. The fact that the function
\begin{equation*}
(y,k)\rightarrow f(x,y,l)s_{x-\rho}(y_i+k)
\end{equation*}
is bounded is elementary.
\end{demo}

Let 
\begin{equation}\label{simplexe2}\index{$A_q$}
A_q= \{ y''=(y_1''\geq\ldots\geq y_n''\geq 0)),\,\, \prod_{i=1}^ny_i''=1,\,\, 
\sum_{i=1}^n y_i''=nq^{-1} \}
\end{equation}
To $y'\in\Sigma$ we associate an element $y''=\phi (y')\in A_q$ defined by
$y_i''=q^{-1}n(y_i'+\alpha )/(1+n\alpha )$ where $\alpha$ is the only
real number such that $\prod_{i=1}^n(y_i'+\alpha )/(1+n\alpha )=1$.
The fact that this map is well defined (i.e. the fact that the
real number $\alpha$ is unique) is a consequence of
the proof of \ref{tcl}, in which it is showed that
$\alpha\rightarrow \prod_{i=1}^n(y_i'+\alpha )/(1+n\alpha )$ is non increasing.
The fact that it is continuous is a consequence of the continuity of the 
roots of a polynomial with respect to its coefficients.
The map $\phi$ admits a left and right
inverse, that, to an element of $y''$ of $A_q$ associates $y'$ defined by 
$y_i'=(y_i''-y_n'')/(\sum_j y_j''-y_n'')$. It is also continuous.
As a summary, we have

\begin{Lem}\label{homeo}
The map $\phi : \Sigma \rightarrow A_q$ is a homeomorphism.
\end{Lem}

The key to the proof of the main result of this section
is a precise understanding of the asymptotics
of the summands of the kernel $u$.

\begin{Lem}\label{tcl}
For all $d$ it is possible to choose two integers
$a_d<b_d$ such that 
\begin{equation*}
\begin{split}
U^{|y^d|}(y^d)\sim \sum_{k=a_d}^{b_d}P^{|y_d|+kn}(y^d)\\
(q/n)^n(|y^d|+kn)\prod_{i=1}^n(y^d_i+k)^{-1}\sim 1\\
\end{split}
\end{equation*}
uniformly in $k\in [a_d,b_d]$ as $d\rightarrow\infty$.
\end{Lem}

\begin{demo}
Let 
\begin{equation*}
f_k(y)=p^{|y|+(k+1)n}(y)/p^{|y|+kn}(y)=(q/n)^n
\frac{(|y|+(k+1)n)!}{(|y|+kn)!\prod_{i=1}^n(y_i+k+1)}
\end{equation*}
This function is defined a priori only for 
$k\in\mathbb{N}^*$, but it admits a natural extension
on the index set $k\in\mathbb{R}_+^*$.

As $y^d\rightarrow y'$, the function
$g_d: t\rightarrow f_{t |y^d|}(y^d)$ 
converges pointwise on $]0,\infty [$ towards 
\begin{equation*}
g_{\infty}: t\rightarrow\frac{q^n(1/n+t )^n}{\prod_{i=1}^n(y_i'+t )}
\end{equation*}
This function is strictly decreasing.
Indeed, its logarithmic derivative is 
\begin{equation*}
\frac{g_{\infty}'}{g_{\infty}} :
t \rightarrow \frac{n}{1/n+t}-\sum_{i=1}^n\frac{1}{y_i'+t}
\end{equation*}
and the inequality between harmonic mean and natural mean implies
that this logarithmic derivative is always $<0$.

Let $[a,b]$ be a closed subinterval of $]0,\infty [$. For $d$ large enough,
$t \rightarrow g_d(t )$ admits a logarithmic derivative 
that is non positive everywhere on $[a,b]$.
Indeed, the function
\begin{equation*}
(t , d)\rightarrow \frac{\partial}{\partial t}\log g_d(t )
\end{equation*}
is easily seen to be a continuous function on the set
$[a,b]\times (\mathbb{N}\cup \{+\infty \})$. 
Therefore for $d$ large enough, $g_d$
is non-increasing.
By Dini's theorem this implies that the convergence 
of $g_d$ towards $g_{\infty}$ holds uniformly on
compact subsets of $]0,\infty [$.

Let $\varepsilon_d=|y_d|^{-1/3}$ and let $[a_d,b_d]$ be the greatest
interval such that $|f_k(y_d)-1|\leq\varepsilon_d$ for all $k\in [a_d,b_d]$.
By the property of uniform convergence on compact subsets and by the 
non increasing property of the limit, this interval is well defined
and there exist non-negative constants $C_1$ and $C_2$ such
that $C_1|y^d|^{2/3}\leq |b_d-a_d|\leq C_2|y^d|^{2/3}$. 
For $k\in [0,a_d-1]$, we have
\begin{equation*}
P^{|y_d|+kn+n}(y^d)/P^{|y_d|+kn}(y^d)\geq 1+|y_d|^{-1/3}
\end{equation*}
and for $k\geq b_d$, 
\begin{equation*}
P^{|y_d|+kn+n}(y^d)/P^{|y_d|+kn}(y^d)\leq 1-|y_d|^{-1/3}
\end{equation*}
An immediate recursion together with a geometric series summation argument
shows that
\begin{equation*}
\sum_{k=0}^{a_d-1}P^{|y_d|+kn}(y^d)\leq P^{|y_d|+a_dn}(y^d)|y_d|^{1/3}
\end{equation*}
and
\begin{equation*}
\sum_{k\geq b_d} P^{|y_d|+kn}(y^d)\leq P^{|y_d|+b_dn}(y^d)|y_d|^{1/3}
\end{equation*}

Let $[a^1_d,b^1_d]$ (resp.$[a^2_d,b^2_d]$ )be the greatest
interval such that $f_k(y_d)\in [1,1+\varepsilon_d]$
(resp.$f_k(y_d)\in [1-\varepsilon_d,1]$).
There is a non-negative constant $C_3$ such that $[a^1_d,b^1_d]$ and
$[a^2_d,b^2_d]$ are of length more than $C_3 |y^d|^{2/3}$.
Furthermore, by the definition of $f_k$ 
there exists an index $i\in\{1,2\}$ such that for any $k\in[a^i_d,b^i_d]$, 
\begin{equation*}
P^{|y_d|+kn}(y^d)\geq\max 
\{P^{|y_d|+a_dn}(y^d),P^{|y_d|+b_dn}(y^d) \}
\end{equation*}
This shows that $U^{|y^d|}(y^d)\sim \sum_{k=a_d}^{b_d}P^{|y_d|+kn}(y^d)$.
\end{demo}

\begin{Prop}\label{martin-estimate}
Let $(y^d)_{d\in\mathbb{N}}$ be a sequence of $\overset{o}{W}$ converging towards
an element $y'\in\Sigma$ as above. 
Then 
\begin{equation*}
\lim_d \frac{\overset{o}{U}(x,y^d)}{\overset{o}{U}(\rho,y^d)}=s_{x-\rho}(y_i'')
\end{equation*}
where $y''=\phi (y')$, as defined in Lemma \ref{homeo}.
\end{Prop}

\begin{demo}
Lemma \ref{lem-schur} and \ref{tcl} imply that 
\begin{equation}\label{a1}
\overset{o}{U}\sim (q/n)^{-|\rho|}\sum_{k\in [a_d,b_d]}P^{|y^d-x|+kn+c_{x,y^d}}(y^d)V(y^d_i)
s_{x-\rho}(y_i'')(|y^d|+kn)^{-|\rho|}
\end{equation}
as $d$ tends to infinity.
Equivalently, we have
\begin{equation}\label{a2}
\overset{o}{U}(\rho,y^d)\sim (q/n)^{-|\rho|}\sum_{k\in [a_d,b_d]}P^{|y^d-\rho|+kn+c_{\rho,y^d}}(y^d)V(y^d_i)
(|y^d|+kn)^{-|\rho|}
\end{equation}
The asymptotics of \ref{a1} (resp. \ref{a2}) do not change if one replaces
the range $[a_d,b_d]$ of the sum by $[a_d-c_{x,y^d}/n,b_d-c_{x,y^d}/n]$
(resp. $[a_d-c_{e,y^d}/n,b_d-c_{e,y^d}/n]$~) 
because $c_{x,y^d}$ takes only a finite number of values (contained in
$[0,|x|]$). Furthermore $k\rightarrow (|y^d|+kn)^{-|e|}$ is a rational fraction,
so Equation \ref{a1} (resp. \ref{a2}) still holds if one replaces 
$c_{x,y^d}$ by $0$ in the r.h.s. (resp. $c_{e,y^d}$ by $0$).
This implies that
\begin{equation*}
\frac{\overset{o}{U}(x,y^d)}{\overset{o}{U}(e,y^d)}\sim s_{x-\rho}(y_i'')
\end{equation*}
\end{demo}

This proves the main part of Theorem \ref{martinbd}, namely that 
the Martin compactification of $W$ is the compactification by $\Sigma$ that we defined in equation
\ref{simplexe1}.
In particular, the minimal harmonic functions are of the kind 
$x\rightarrow s_{x-\rho}(y_i'')$, where $y_i''\in A_q$.

\subsection{Extremal Martin Boundary}\label{mex}

This part is devoted to computing $MB^{min}$, and showing that
$MB^{min}=MB$.
Let $\mathfrak{sl}_n$ be the complex Lie algebra of $SL(n)$. 

$SU(n)$ admits a natural left action by conjugation on 
$\widehat{M(SU(n))}$, which we denote by $\Ad$.
The normalized trace $\tr$ of the fundamental representation.
extends by linearity and continuity 
to a positive linear functional on $\widehat{M(SU(n))}$.
This allows to define the positive convolution operator
\begin{equation}
\widetilde{P}=(\tr \otimes id)\D
\end{equation}

Let $\widehat{Z(SU(n))}$ be the center of $\widehat{M(SU(n))}$. Since $\tr$
commutes with the action $\Ad$, so does $\widetilde{P}$.  
The algebra
$\widehat{Z(SU(n))}$ is the fixed point algebra of $\Ad$, therefore
$\widetilde{P}$  leaves $\widehat{Z(SU(n))}$ invariant
and defines a new submarkovian operator on $\overset{o}{W}_n$ (upon the obvious
identification of $\mathbb{C}^{\overset{o}{W}_n}$ with 
$\widehat{Z(SU(n))}$).

\begin{Prop}\label{min-th}
The minimal harmonic functions with respect to the operator $P$ 
(resp. $\widetilde{P}$) on $\overset{o}{W}$ are the functions
\begin{equation*}
x\rightarrow s_{x-\rho}(y_1'',\ldots ,y_n'')
\end{equation*}
(resp. $x\rightarrow s_{x-\rho}(y_1'',\ldots ,y_n'')/s_{x-\rho}(1)$),
where $(y_1'',\ldots ,y_n'')$ run in $A_q$.
\end{Prop}

\begin{demo}
It is equivalent to have
$f : \overset{o}{W} \rightarrow \mathbb{R}_+^*$ harmonic with
respect to $\overset{o}{P}$, and 
$\tilde{f} : x\in\overset{o}{W} \rightarrow f(x)/s_{x-\rho}(1)$ 
harmonic with respect to $\widetilde{P}$.

By \cite{MR94a:46091}, for any positive harmonic element $f$,
there exists a finite positive
measure $\mu_f$ on the set of positive elements of $SL(n)$
(we call it $SL(n)_+$),
such that $f=\int_{SL(n)}i(A)d\mu (A)$. 
Since $f$ is invariant under $\Ad$, $\mu_f$ is also invariant under
$\Ad$.
Let $\widetilde{\mu_f}$ be the image of $\mu_f$ under the canonical
projection of $SL(n)_+$  onto $SL(n)_+/\Ad$. 
The quotient space $SL(n)_+/\Ad$ contains naturally $A_q$, and the support
of $\widetilde{\mu_f}$ is a subset of $A_q$ upon this inclusion.

Conversely, a finite positive measure $d\mu$ in $A_q$ represents an element
of $\widehat{Z(SU(n))}$. Indeed, let $Y$ be the matrix $diag (y_1,\ldots ,y_n)$.
Recall that $i$ was defined at equation \ref{def-i}.
Then,
\begin{equation*}
\int_{SU(n)}\int_{A_q} i(UYU^*)dUd\mu 
\end{equation*}
defines an harmonic element in $\widehat{Z(SU(n))}$.
This defines a one to one correspondence between 
positive finite measures on
$A_q$ and harmonic elements of $\widehat{Z(SU(n))}$.

Let
$y=(y_1\geq\ldots\geq y_n)\in A_q$. 
Then for any $U\in SU(n)$ and for any $x\in\overset{o}{W}$,
we have $s_{x-\rho}(UYU^*)=s_{x-\rho}(y_i)$. Therefore, denoting by $p_x$ the minimal
central idempotent of $\mathfrak{U(sl_n)}$, the element
$\int_{SU(n)} i(UMU^*)dU $
is harmonic, central and by equality of traces,
\begin{equation*}
p_x\int_{SU(n)} i(UMU^*)dU=p_x s_{x-\rho}(y)/s_{x-\rho}(1)
\end{equation*}
Therefore, to each element $y$ of $A_q$ corresponds the 
harmonic function
\begin{equation*}
x\rightarrow s_{x-\rho}(y)/s_{x-\rho}(1)
\end{equation*}

This proves $f$ has to be
$x\rightarrow s_{x-\rho}(y)/s_{x-\rho}(1)$ for some $y\in A_q$. 
\end{demo}

\section{Quantum Random Walks on $\widehat{SU_q(n)}$ and \\ Abelian
subalgebras}\label{fusion}

In this part, we apply the results of section \ref{nsfa} to a quantum
random walk on the dual of $SU_q(n)$ that was first introduced in
\cite{Iz}.

  \subsection{Quantum compact groups}\label{defqgp}

We start with the definition of matrix pseudogroup due to Woronowicz (see
\cite{MR90e:22033} and \cite{MR88m:46079}).
Let $A$ be a $C^*$ -algebra with unit. The set $\mathbb{M}_N(A)$ of matrices with
entries belonging to $A$ is identified with the $C^*$-algebra $B(\mathbb{C}^N)\otimes A$.
A pair  $(u=(u_{ij})\in\mathbb{M}_N(A), A)$ is said to be a 
{\it compact matrix pseudogroup \it} iff \index{groupe quantique|(}
\begin{itemize}
\item the $*$-subalgebra $\mathcal{A}$ generated by matrix elements of $u$ is
dense in $A$.\index{$\mathcal{A}$}
\item there exists a $C^*$-homomorphism 
\begin{equation*}
\Delta : A\rightarrow A\otimes A
\end{equation*}
such that
\begin{equation*}
\Delta (u_{ij})=\sum_{k=1}^n u_{ik}\otimes u_{kj} 
\end{equation*}
\item
$u$ is invertible and there exists a linear antimultiplicative mapping 
\begin{equation*}
\kappa : \mathcal{A}\rightarrow\mathcal{A}
\end{equation*}
such that $\kappa (\kappa (a^*)^*)=a$ and
\begin{equation*}
({\rm id \rm} \otimes \kappa) u=u^{-1}
\end{equation*}
\end{itemize} 

An element $w=(w_{ij})\in\M{n}\otimes A$ is called a
{\it unitary corepresentation\it} if the following holds:
\begin{equation*}
\Delta w_{ij}=\sum_k w_{ik}\otimes w_{kj}
\end{equation*}
A vector space $V$ with basis $v_i$ and with a map
$\Phi : V\rightarrow V\otimes A$ is called a comodule
if there exists a corepresentation of $A$ such that
\begin{equation*}
\Phi v_j=\sum_k v_k\otimes w_{kj}
\end{equation*}
For example, $vect \{w_{i1},\ldots ,w_{ik}\}$ is a comodule if
$w=(w_{ij})$ is a unitary corepresentation.
It is possible to define in an obvious way a notion of subcomodule,
irreducible comodule, and equivalent comodules.

In this paper we shall focus on the specific example of
$A(SU_q(n))$. It is the universal $C^*$ -algebra generated by $n^2$ elements
$u_{kl}$ ($k,l=1,2,\ldots ,n$) such that \index{$SU_q(n)$}
\begin{itemize}
\item
\begin{equation}\label{defsuq}
\sum_k u_{kl}^*u_{km}=\delta_{lm}I \,\, , \,\,
\sum_k u_{mk}u_{lk}^*=\delta_{lm}I
\end{equation}
\item
\begin{equation}
\sum_{k_1,\ldots ,k_n}u_{l_1k_1}\ldots u_{l_nk_n}E_{k_1,k_2,\ldots ,k_n}=
E_{l_1,l_2,\ldots ,l_n}I
\end{equation}
\end{itemize}
where, for $q\in ]0,1]$, 
\begin{equation}
E_{i_1,i_2,\ldots ,i_n}=
\begin{cases}
0 {\rm \,\, if \,\, i_k=i_l \,\, for \,\, some \,\, k\neq l \rm}\\
(-q)^{I(i_1,i_2,\ldots ,i_n)} \,\, {\rm otherwise \rm}
\end{cases}
\end{equation}
with $I(i_1,i_2,\ldots ,i_n)$ denoting the number of inversed pairs in the
sequence $(i_1,i_2,\ldots ,i_n)$.
Then $SU_{q}(n)=(A(SU_q(n)),u)$ is a compact matrix pseudogroup. 
Furthermore, for $q = 1$, it coincides with the algebra of continuous functions on the
classical $SU(n)$ group.

For any matrix pseudogroup,
there exists a unique state $h$ called {\it Haar measure\it}, satisfying 
\begin{equation*}
h(x)\cdot 1= (h\otimes id)\cdot \Delta (x) = (id \otimes h) \cdot \Delta (x), x\in A
\end{equation*}

The state $h$ is faithful in the case of $SU_q(n)$.
Let $(\pi_h, H_h, \Omega_h)$ be the GNS triple of $h$, and $\Lambda_h$ the 
natural map from $A(SU_q(n))$ to $B(H_h)$. The {\it multiplicative unitary \it}
is defined as the bounded extension of the following operator:
\begin{equation*}\index{unitaire multiplicatif}
V(\Lambda_h(x)\otimes \xi)=\Delta (x) (\Omega_h\otimes\xi ), x\in A, \xi\in H_h
\end{equation*}
$V$ is unitary and
satisfies the following pentagon equation: (see \cite{MR94e:46127})
\begin{equation*}
V_{12}V_{13}V_{23}=V_{23}V_{12}
\end{equation*}
The dual von Neumann algebra $M(SU_q(n))$ is the bicommutant in $B(H_h)$ of the
set $\{ (id \otimes \mu) V \}$ where $\mu$ runs over $B(H_h)_*$.
The dual coproduct is defined by
\begin{equation}\index{$\D$}
\D(x)=V^* (1\otimes x)V
\end{equation}
therefore $M(SU_q(n))$ is endowed with a Hopf-von Neumann algebra structure.
This von Neuman algebra is well understood: since
the representation theory of $SU_q(n)$ is the same as
that of $SU(n)$, it has the same von Neumann algebra
structure as the von Neumann algebra of $SU(n)$, therefore is isomorphic to
\begin{equation*}
\bigoplus_{x\in W}\M{d_x}
\end{equation*}
\index{groupe quantique|)}

For any representation $s\in W$, let $\{f_z\}_{z\in\mathbb{C}}$ be the family of 
Woronowicz characters. For its definition and basic properties, we refer to 
\cite{MR88m:46079}.\index{caractères!de Woronowicz}
We only need to know that there exists a unique positive $\rho\in\M{n}$ with
normalized trace, such that $f_z=\rho^z$, and $\rho$ intertwines the fundamental
representation with its double contragredient.
Let $s$ be the fundamental representation, $\mu =\tr_s(\rho \cdot)$, and
\begin{equation}\label{defq}
\tilde{P}:=(id\otimes \mu )\D
\end{equation}
The operator $\tilde{P}$ is completely positive. It has already been considered
by Izumi (see \cite{MR1697285}, \cite{MR99g:46093} and his preprint \cite{Iz}). 
It leaves invariant the center \index{$Z(A)$}
$Z (SU_q(n))$ of $M(SU_q(n))$ (this is a consequence
of \cite{Iz}, Lemma 3.2, (3) ).
In order to prove this, one needs to define the morphism $\Phi$:\index{$\Phi$}
\begin{equation*}
\begin{split}
M(SU_q(n) )\rightarrow M(SU_q(n) )\otimes L^{\infty} (SU_q(n) )\\
\Phi (x)=U (x\otimes 1) U^*
\end{split}
\end{equation*}
This is a von Neumann action of $SU_q(n)$ on $M(SU_q(n) )$.
The center $Z(SU_q(n) )$ is invariant under this action, therefore
it makes sense to restrict it to Peter-Weyl blocks and to extend it
to $\widehat{M(SU_q(n) )}$. 

$\tilde{P}$ intertwines $\Phi$, therefore it leaves invariant the center
$Z(SU_q(n))$ of $M(SU_q(n))$. We define the twisted integer 
\begin{equation*}\index{entier déformé}\index{$[n]_q$}
[n]_q=\frac{q^n-q^{-n}}{q-q^{-1}}
\end{equation*}
With this, one can show

\begin{Prop}
In the canonical basis of $\M{n}$, we have
\begin{equation}
\rho = \frac{1}{[n]_q}
\left(
\begin{matrix}
q^{-n+1} & 0 & \ldots & 0 \\
0 & q^{-n+3} & \ddots & \vdots \\
\vdots & \ddots & \ddots & 0 \\
0 & \ldots & 0 & q^{n-1}
\end{matrix}
\right)
\end{equation}
\end{Prop}

It is a standard computation. See for example \cite{MR88m:46079}.

\subsection{Restrictions and $h$-processes}

Let $A(T^{n-1})$ be the $C^*$-algebra of continuous function on the torus
$T^{n-1}$. It is the universal Abelian $C^*$-algebra generated by
the $n$ unitaries $u_1,\ldots ,u_n$ satisfying $u_1\ldots u_n=1$. 
The coproduct $\D u_i=u_i\otimes u_i$ defines a $C^*$ -Hopf algebra structure.

Let $\phi$ be the algebra morphism $A(SU_q(n))\rightarrow A(T^{n-1})$
such that $\phi (u_{ij})=\delta_{ij}u_i$
with the induced relations. Obviously $\phi$ is a morphism of 
$C^*$ -Hopf algebras.
This allows to define a subalgebra 
$M(T^{n-1})$ of $M(SU_q(n))$. 
It is the von Neumann algebra generated by 
$\{ (id\otimes \mu)V \}$, where $\mu$ runs over the characters of $A(T^{n-1})$
seen as elements of $B(H_h)_*$.
One can show that it is isomorphic to 
the group von Neumann algebra of $T^{n-1}$ and has a Hopf-von Neumann 
algebra structure. 
Therefore it makes sense to restrict $P$ defined at Equation
\ref{defq} to this
von Neumann subalgebra. 

Furthermore, $\tilde{P}$ is invariant under the adjoint action $\Phi$, 
thus it leaves invariant the center $Z(SU_q(n))$.\index{$\Phi$}

Therefore, $\tilde{P}$ restricts to two natural Abelian subalgebras of $M(SU_q(n))$.
$M(T^{n-1})$ is in natural correspondance with $L^{\infty}(L)$ defined in
section \ref{nsfa}
and $Z(SU_q(n))$ with $L^{\infty}(\overset{o}{W_n})$. Thus $P$ induces a classical 
Markov chain on $L$ and an other one on $\overset{o}{W_n}$. It is of natural interest
to investigate a (probabilistic) link between these two Markov chains.

This has already been done by Biane in \cite{MR93a:46119} in the case $q=1$.
He shows that under the embedding of $\overset{o}{W_n}$ into $L$ described in the preceeding
section \ref{nsfa}, 
the Markov chain on $\overset{o}{W_n}$ is obtained from that on $L$ by conditioning
in Doob's sense the Chain on $L$ not to leave $\overset{o}{W_n}$. The chain on $L$
itself is a nearest neighbor centered random walk.

Consider, for $q<1$, the point of $\Sigma\subset S^{n-2}$
\begin{equation}
y_q=\sum_{i=1}^n q^{n+1-2i}e_i/||\sum_{i=1}^n q^{n+1-2i}e_i||\in S^{n-2}
\end{equation}

\begin{Prop}
The random walk induced by the restriction of $\tilde{P}$ to $M(T^{n-1})$
corresponds on the lattice $L$, to
the convolution by the probability measure 
\begin{equation}
\mu (e_i) =\frac{q^{-n+2i-1}}{[n]_q}
\end{equation}
It is obtained by conditioning, in Doob's sense, the submarkov random walk
\begin{equation}
\tilde{\mu} (e_i) =\frac{1}{[n]_q} 
\end{equation}
not to die and to converging towards the point $y_q$.
\end{Prop}

\begin{demo}
It is enough to remark that the function 
\begin{equation*}
h : \begin{cases} L\rightarrow \mathbb{R}\\
\sum x_ie_i\rightarrow q^{\sum_i x_i(-n+2i-1)}
\end{cases}
\end{equation*}
is well defined and harmonic with respect to $\tilde{\mu}$.
Furthermore the Doob conditioning of the convolution operator by $\tilde{\mu}$ 
with respect to $h$ is the convolution operator by $\mu$ and by the
law of large numbers, this random walk almost surely tends to $y_q$.
\end{demo}

Now we focus on the restriction to the center. 

\begin{Prop}\label{transpro}
The transition probability of the random walk restricted to the center
is 
\begin{equation}
p_{l,l'}=
\begin{cases} \frac{s_{l'}(q^{-n+1}, q^{-n+3},\ldots , 
q^{n-1})}{s_{l}(q^{-n+1}, q^{-n+3},\ldots , q^{n-1})[q]_n} 
{\rm \,\, if \,\,\exists i,\,\, l'=l+e_i \,\, and \,\,l'\in\overset{o}{W}  \rm}\\
0 {\rm \,\,otherwise \rm}
\end{cases}
\end{equation}
\end{Prop}

\begin{demo}
This is a consequence of  
\cite{Iz} and \cite{MR99h:46116}
together with the fact that the quantum dimension of
the representation $l$ is 
\begin{equation*}
s_l(q^{-n+1}, q^{-n+3},\ldots , q^{n-1})
\end{equation*}
This last fact results from the fact that the representation
theory of $SU_q(n)$ and $SU(n)$ are the same, that the maximal
torus remains non-deformed, and that we therefore have a
Weyl character formula. 
\end{demo}

\begin{Thm}\label{transpro-b}
The restriction of $\widetilde{P}$ to the center corresponds to
the random walk on $L$ with increment $\alpha \sum \delta_{e_i}/n$
with $\alpha = n/[n]_q$, conditioned to converging towards
the point $y_q$, not to die, and not to hit $\partial W$.
\end{Thm}

\begin{demo}
According to the previous section, Proposition \ref{min-th} and Theorem \ref{martinbd},
$s\rightarrow s_x(q^{-n+1},q^{-n+3},\ldots ,q^{n-1})$ is harmonic
with respect to $\overset{o}{P}$ and in the Martin theory, it corresponds
to the point $y_q$. This implies that the Doob conditioning of  
$\overset{o}{P}$ with respect to
this harmonic function corresponds to conditioning to tend towards $y_q$
and remain inside $\overset{o}{P}$. 
But this operator is also $\tilde{P}$, as it was computed in Proposition \ref{transpro}.
\end{demo}

\bibliographystyle{alpha}
\bibliography{article9092-1}

\end{document}